\documentstyle[12pt]{article}
\newtheorem{theorem}{Theorem}

\newtheorem{proposition}{Proposition}
\newtheorem{definition}{Definition}
\newtheorem{corollary}{Corollary}

\begin{document}

\title{Cohomology $C_{\infty}$-algebra and
Rational Homotopy Type}
\date{}
\author{Tornike Kadeishvili}

\maketitle


\abstract{In the rational cohomology of a 1-connected space a
structure of $C_{\infty}$-algebra is constructed and it is shown
that this object determines the rational homotopy type}

\footnotetext[1]{The research described in this publication was
made possible in part by Award No. GNSF/ST06/3-007 of the Georgian
National Science Foundation}



\section{INTRODUCTION}

Usually invariants of algebraic topology are not {\it complete}:
the isomorphism of invariants does not guarantee the equivalence
of spaces. The invariants which carry richer algebraic structure
contain more information about the space. For example the
invariant "cohomology algebra" allows to distinguish spaces, which
can not be distinguished by the invariant "cohomology groups".

Let us assume that all $R$-modules $H^*(X,R)$ are free. In [24,25]
we obtain an $A_\infty$-algebra structure on $H^*(X,R)$. This
structure consists of a collection of operations
$$
\{m_i: H^*(X,R)\otimes ... (i\ times)...\otimes H^*(X,R)\to
H^*(X,R),\ i=2,3,...\}.
$$
In fact this structure extends the usual structure of cohomology
algebra: the first operation $m_2:H^*(X,R)\otimes H^*(X,R)\to
H^*(X,R)$ coincides with the cohomology multiplication.

The cohomology algebra equipped with this additional structure,
which we call cohomology $A_{\infty}$-algebra, carries more
information about the space, than the cohomology algebra. For
example just the cohomology algebra $H^*(X,R)$ does not determine
cohomology of the loop space $H^*(\Omega X,R)$, but the cohomology
$A_{\infty}$-algebra $(H^*(X,R),\{m_i\})$ does. Dually, the
Pontriagin ring $H_*(G)$ does not determine homology $H_*(B_G)$ of
the classifying space, but the homology $A_{\infty}$-algebra
$(H_*(G),\{m_i\})$ does.

These $A_{\infty}$-algebras has several applications in the
cohomology theory of fibre bundles too, see \cite{Kad80}.

But this invariant also is not complete. One can not expect the
existence of more or less simply complete algebraic invariant in
general case but for the rational homotopy category there are
various complete homotopy invariants (algebraic models):

(i) The model of Quillen \cite{Qui69} $L_X$, which is a
differential graded Lie algebra;

(ii) The minimal model of Sullivan \cite{Sul75} $M_X$, which is a
commutative graded differential algebra;

(iii) The filtered model of Halperin and Stasheff \cite{HS}
$\Lambda X$, which is a filtered commutative graded differential
algebra.

The rational cohomology algebra $H^*(X,Q)$ is not a complete
invariant even for rational spaces: two spaces might have
isomorphic cohomology algebras, but different rational homotopy
types.

The main result of this paper is the construction of a complete
rational homotopy invariant: {\it the cohomology
$C_{\infty}$-algebra}.

This notion of $C_{\infty}$-algebra is the commutative version of
the Stasheff's notion of $A_{\infty}$-algebra. It was mentioned in
\cite{Smi80}; in \cite{Kad88} it was called commutative
$A_{\infty}$-algebra and was denoted as $CA_{\infty}$; in
\cite{Markl} it was called {\it balanced} $A_{\infty}$-algebra;
the modern notation $C_{\infty}$-algebra was introduced in
\cite{GJ94}.

We show that in the rational case on cohomology $H^*(X,Q)$ arises
a structure of $C_{\infty}$-algebra $(H^*(X,Q),\{m_i\})$. The main
application of this structure is following: it completely
determines the rational homotopy type, that is 1-connected spaces
$X$ and $X'$ have the same rational homotopy type if and only if
their cohomology $C_{\infty}$-algebras $(H^*(X,Q),\{m_i\})$ and
$(H^*(X',Q),\{m'_i\})$ are isomorphic.

We present also several applications of this complete rational
homotopy invariant to some problems of rational homotopy theory.

The $C_\infty$-algebra structure in homology of a
\emph{commutative} dg algebra and the applications of this
structure in rational homotopy theory was actually presented in
hardly available small book \cite{Kad93} (see also the preprint
\cite{Kad88h}).

Applications of cohomology $C_\infty$-algebra in rational homotopy
theory are inspired by the existence of Sullivan's commutative
cochains $A(X)$ in this case. The cohomology $C_\infty$-algebra
$(H^*(X,Q),\{m_i\})$ carries the same amount of information as
$A(X)$ does. Actually these two objects are equivalent in the
category of $C_\infty$-algebras.

Outside of rational category generally we do not have commutative
cochains, so some additional structures, such as Steenrod
$\smile_i$ products, and much more, must be involved. For example
as the first step one should add the operations which form so
called homotopy G-algebra structure (in fact the little square
operad) (\cite{GV95}, \cite{Kad07}). These in fact are cochain
operations which control interaction between $\smile$ and
$\smile_1$ products. Next, some new operations which control
interaction between $\smile$ and $\smile_i,\ i=1,2,3,...$ products
show up (\cite{Kad03}). Next must be operations which control
interaction between $\smile_i$ and $\smile_j$ products, etc.

We presume that finally we obtain some specific $E_\infty$ algebra
structure on singular cochains, see \cite{Hinich}, \cite{Mac},
\cite{Fresse}.

The final achievement in this direction is Mandel's result: the
$E_\infty$-algebra structure on cochain algebra determines (in
some cases) the homotopy type.

In rational case $E_\infty$ operad can be replaced by commutative
operad $\cal{C}$ acting on appropriate cochains. And in order to
step from cochains to cohomology we replace $\cal{C}$ be the
operad $C_\infty$.

\section{$A_{\infty}$-algebras}

The notion of $A_{\infty}$-algebra was introduces by J. Stasheff
\cite{Sta63}. This notion generalizes the notion of differential
graded algebra (dga).

\begin{definition}\label{ainf}
An $A_{\infty}$-algebra is a graded module $M=\{M^k\}_{k\in Z}$
equipped with a sequence of operations
$$
\{m_i:M\otimes ...(i-times)...\otimes M\to M, i=1,2,3,. . . \}
$$
satisfying the conditions $ m_i((\otimes^iM)^q)\subset M^{q-i+2}
$, that is $deg\ m_i=2-i$, and
\begin{equation}
\label{ainfalg}
\begin{array}{l}
\sum_{k=0}^{i-1} \sum_{j=1}^{i-k} \pm \\
m_{i-j+1}(a_1\otimes . . . \otimes a_k\otimes  m_j(a_{k+1}\otimes
. . . \otimes a_{k+j})\otimes . . . \otimes a_i)=0 .
\end{array}
\end{equation}
\end{definition}
In fact for an $A_{\infty}$-algebra $(M,\{m_i\})$ first two
operations form a {\it nonassociative} dga $(M,m_1,m_2)$ with
differential $m_1$ and multiplication $m_2$ which is associative
just up to homotopy and the suitable homotopy is the operation
$m_3$.

\begin{definition}
 A  morphism of $A_{\infty}$-algebras
$$ \{f_i\}:
(M,\{m_i\})\to (M',\{m'_i\})
$$
is a sequence
$
 \{f_i:\otimes^iM\to M', i=1,2,. . . ,\ deg\ f_1=1-i \}
$ such that
 \begin{equation}
\label{morphism}
\begin{array}{l}
\sum_{k=0}^{i-1} \sum_{j=1}^{i-k} \pm \\
f_{i-j+1}(a_1\otimes . . . \otimes a_k \otimes m_j(a_{k+1}\otimes
. . . \otimes a_{k+j})\otimes . . .  \otimes a_i)=\\
\sum_{t=1}^{i} \sum_{k_1+...+k_t=i}\pm \\
m'_t(f_{k_1}(a_1\otimes. . . \otimes a_{k_1})\otimes. . . \otimes
f_{k_t}(a_{i-k_t+1}\otimes. . . \otimes a_{i})).
\end{array}
\end{equation}
\end{definition}
The composition of $A_{\infty}$ morphisms
$$
\{h_i\}:(M,\{m_i\})\stackrel{\{f_i\}}{\longrightarrow}
(M',\{m'_i\})\stackrel{\{g_i\}}{\longrightarrow} (M'',\{m''_i\})
$$
is defined as
\begin{equation}
\label{composition}
\begin{array}{l}
 h_n(a_1\otimes...\otimes
a_n)=\sum_{t=1}^{n}\sum_{k_1+...+k_t=n}\\
g_{n}(f_{k_1}(a_{1}\otimes...\otimes a_{k_1})\otimes ...\otimes
f_{k_t}(a_{n-k_t+1}\otimes...\otimes a_{n}).
\end{array}
\end{equation}
The bar construction argument (see (\ref{bartilda}) bellow) allows
to show that so defined composition satisfies the
condition(\ref{morphism}).

For a morphism $\{f_i\}:(M,\{m_i\})\to (M',\{m'_i\})$ the first
component $f_1:(M,m_1)\to (M',m'_1)$ is a chain map which is {\it
multiplicative} just up to homotopy and the suitable homotopy is
the map $f_2$.

$A_{\infty}$ algebra of type $(M,\{m_1,m_2,0,0,...\})$ is a dga
with the differential $m_1$ and strictly associative
multiplication $m_2$. Furthermore, a morphism of such
$A_{\infty}$-algebras of type $\{f_1,0,0,...\}$ is a strictly
multiplicative chain map. Thus the category of dg algebras is the
subcategory of the category of $A_{\infty}$-algebras.


\section{$C_{\infty}$-algebras}


The {\it shuffle} product $ \mu_{sh}:M^{\otimes m}\otimes
M^{\otimes n}\to M^{\otimes (m+n)} $ is defined as
\begin{equation}
\label{shuffle}
 \mu ((a_1\otimes ... \otimes a_n)\otimes
(a_{n+1}\otimes ... \otimes a_{n+m}))= \sum \pm a_{\sigma
(1)}\otimes ... \otimes a_{\sigma (n+m)},
\end{equation}
where summation is taken over all $(m,n)$-shuffles, that is over
all permutations of the set $(1,2,...,n+m)$ which satisfy the
condition: $i<j$ if $1\leq \sigma (i)<\sigma (j)\leq n$ or
$n+1\leq \sigma (i)<\sigma (j)\leq n+m$.

\begin{definition} (\cite{Smi80}, \cite{Kad88},\cite{Markl}, \cite{GJ94}) A $C_{\infty}$-algebra is an
$A_{\infty}$-algebra $(M,\{m_i\})$ which additionally satisfies
the following condition: each operation $m_i$ disappears on
shuffles, that is for $a_1,...,a_i\in M$ and $k=1,2,...,i-1$
\begin{equation}
\label{m0onshuffle} m_i(\mu_{sh}((a_1\otimes ...\otimes
 a_k)\otimes(a_{k+1}\otimes ...\otimes a_i)))=0.
\end{equation}
\end{definition}

\begin{definition}
A morphism of $C_{\infty}$-algebras is defined as a morphism of
$A_{\infty}$-algebras $ \{f_i\}: (M,\{m_i\})\to (M',\{m'_i\}) $
whose components $f_i$ disappear on shuffles, that is
\begin{equation}
\label{f0onshuffle} f_i((\mu_{sh}(a_1\otimes ...\otimes
 a_k)\otimes(a_{k+1}\otimes ...\otimes a_i)))=0.
\end{equation}
\end{definition}
The composition is defined as in $A_{\infty}$ case and the bar
construction argument (see  (\ref{bartilda}) bellow) allows to
show that the composition is a $C_{\infty}$ morphism.

In particular for the operation $m_2$ we have $m_2(a\otimes b\pm
b\otimes a)=0$, so a $C_{\infty}$-algebra of type
$(M,\{m_1,m_2,0,0,...\})$ is a commutative dg algebra (cdga) with
the differential $m_1$ and strictly associative and commutative
multiplication $m_2$. Thus the category of cdg algebras is the
subcategory of the category of $C_{\infty}$-algebras.

\section{Tensor coalgebra}

The notions of $A_{\infty}$ and $C_{\infty}$ algebras can be
interpreted in terms of \emph{differentials} on the tensor
coalgebra.

The {\it tensor coalgebra} of a graded module $V$ is defined as
$$
T^c(V)=R\oplus V\oplus V\otimes V\oplus V\otimes V\otimes V\oplus
...=\sum_{i=0}^\infty V^{\otimes i}
$$
with the comultiplication $\Delta:T^c(V)\to T^c(V)\otimes T^c(V)$
given by
$$
\Delta (a_1\otimes. . . \otimes a_n)=\sum_{i=0}^{n}(a_1\otimes. .
. \otimes a_i)\otimes (a_{i+1}\otimes. . . \otimes a_n).
$$
Tensor coalgebra is the cofree object in the category of graded
coalgebras: for a map of graded modules $\alpha :C\to V$ there
exists unique morphism of graded coalgebras $f_{\alpha}:C\to
T^c(V)$ such that $p_1f_{\alpha}=\alpha$, here $p_n:T^c(V)\to
V^{\otimes n}$ is the clear projection. The coalgebra map
$f_\alpha$ is defied as $ f_\alpha=\sum_k
(\alpha\otimes...\alpha)\Delta^k, $ where $\Delta^k:C\to
C^{\otimes k}$ is the $k$-th iteration of the comultiplication
$\Delta :C\to C\otimes C$, i.e. $\Delta^1=id,\ \Delta^2=\Delta,\
\Delta^k=(\Delta^{k-1}\otimes id)\Delta$.

Tensor coalgebra has similar universal property also for
coderivations, i.e. maps $\partial:C\to C'$ satisfying $\Delta
\partial=(\partial\otimes id+id\otimes\partial)\Delta$. Namely, for
each homomorphism $\beta :T^c(V)\to V$ there exists unique
coderivation $\partial_{\beta}:T^c(V)\to T^c(V)$ such that
$p_1\partial_{\beta}=\beta$. The coderivation $\partial_\beta$ is
defied as $
\partial_\beta=\sum_{k,i} (id\otimes\beta\otimes id)\Delta^3.
$


The {\it shuffle} multiplication $\mu_{sh}:T^c(V)\otimes T^c(V)\to
T^c(V)$, introduced by Eilenberg and MacLane \cite{EM53}, turns
$(T^c(V),\Delta,\mu_{sh})$ into a graded bialgebra.

This multiplication is defined as a graded coalgebra map induced
by the universal property of $T^c(V)$ by $\alpha:T^c(V)\otimes
T^c(V)\to V$ given by $\alpha(v\otimes 1)=\alpha(1\otimes v)=v$
and $\alpha=0$ otherwise. This multiplication is associative and
in fact is given by
$$
\mu_{sh}([a_1,...,a_m]\otimes[a_{i+1},...,a_{n}])=\sum \pm
[a_{\sigma(1)},...,a_{\sigma(n)}]),
$$
where the summation is taken over all $(m,n)$-shuffles.



\subsection{Bar construction of an $A_{\infty}$-algebra}
\label{bartilda}

 Let  $(M,\{m_i\})$ be an $A_{\infty}$-algebra. We consider the
 tensor coalgebra $T^c(s^{-1}M)$ where $s^{-1}M$ is the {\it desuspension}
 of $M$, i.e. $(s^{-1}M)^n=M^{n+1}$. We use the standard notation
$s^{-1}a_1\otimes ... \otimes s^{-1}a_n=[a_1,...,a_n]$. The
structure maps $m_i$ define the map  $\beta:T^c(s^{-1}M)\to
s^{-1}M$ by $\beta [a_1,...,a_n]=[s^{-1}m_n(a_1\otimes ...\otimes
a_n)]$. Extending this $\beta$ as a coderivation we obtain
$d_\beta:T^c(s^{-1}M)\to T^c(s^{-1}M)$ which in fact looks as
$$
d_\beta[a_1,...,a_n]=\sum_k\pm [a_1,...,a_k,m_j(a_{k+1}\otimes
...\otimes a_{k+j}),a_{k+j+1},... a_n].
$$
The defining condition (\ref{ainfalg}) of $A_{\infty}$-algebra
guarantees that $d_\beta d_\beta = 0$. The obtained dg coalgebra
$(T^c(s^{-1}M),d_\beta,\Delta)$ is called {\it bar construction}
of $A_{\infty}$-algebra $(M,\{m_i\})$ and is denoted by
$\tilde{B}(M)$.

For an $A_{\infty}$-algebra of type $(M,\{m_1,m_2,0,0,...\})$ this
bar construction coincides with the ordinary bar construction of
this dga.

A morphism of $A_{\infty}$-algebras $ \{f_i\}: (M,\{m_i\})\to
(M',\{m'_i\}) $ defines a dg coalgebra map of bar constructions
$F=\tilde{B}(\{f_i\})$ as follows: $\{f_i\}$ defines the map
$\alpha:T^c(s^{-1}M)\to s^{-1}M$ by $ \alpha
[a_1,...,a_n]=[s^{-1}f_n(a_1\otimes ...\otimes a_n)]. $ Extending
this $\alpha$ as a coalgebra map we obtain $F:T^c(s^{-1}M)\to
T^c(s^{-1}M)$ which in fact looks as
$$
F[a_1,...,a_n]=\sum \pm [f_{k_1}(a_{1}\otimes...\otimes a_{k_1}),
..., f_{k_t}(a_{n-k_t+1}\otimes...\otimes a_{n})].
$$
The defining condition (\ref{morphism}) of $A_{\infty}$ morphism
guarantees that $F$ is a chain map.

Now we are able to show that the composition of $A_{\infty}$
morphisms is correctly defined: to the composition of morphisms
(\ref{composition}) corresponds the composition of dg coalgebra
maps
$$\tilde{B}((M,\{m_i\}))\stackrel{\tilde{B}(\{f_i\})}{\longrightarrow}
\tilde{B}((M',\{m'_i\}))\stackrel{\tilde{B}(\{g_i\})}{\longrightarrow}
\tilde{B}((M'',\{m''_i\}))$$ which is a dg coalgebra map, thus for
the projection $p_1\tilde{B}(\{g_i\})\tilde{B}(\{f_i\})$, i.e. for
the collection $\{h_i\}$, the condition (\ref{morphism}) is
satisfied.

\subsection{Bar construction of a $C_{\infty}$-algebra}
\label{bartildac}

The notion of $C_{\infty}$-algebra is motivated by the following
observation. If a dg algebra $(A,d,\mu)$ is graded commutative
then the differential of the bar construction $BA$ is not only a
coderivation but also a derivation with respect to the shuffle
product, so the bar construction $(BA,d_\beta,\Delta,\mu_{sh})$ of
a cdga is a dg bialgebra.

By definition the bar construction of an $A_{\infty}$-algebra
$(M,\{m_i\})$ is a dg coalgebra
$\tilde{B}(M)=(T^c(s^{-1}M),d_\beta,\Delta)$.

But if $(M,\{m_i\})$ is a $C_{\infty}$-algebra, then
$\tilde{B}(M)$ becomes a dg bialgebra:
\begin{proposition}
For an $A_{\infty}$-algebra $(M,\{m_i\})$ the differential of the
bar construction $d_\beta$ is a derivation with respect to the
shuffle product if and only if each operation $m_i$ disappears on
shuffles, that is $(M,\{m_i\})$ is a $C_{\infty}$-algebra.
\end{proposition}
\noindent {\bf Proof.} The map $ \Phi:T^c(s^{-1}M)\otimes
T^c(s^{-1}M)\to T^c(s^{-1}M) $ defined as
$\Phi=d_\beta\mu_{sh}-\mu_{sh}(d_\beta\otimes id+id\otimes
d_\beta)$ is a coderivation. Thus, according to universal property
of $T^c(s^{-1}M)$ the map $\Phi$ is trivial if and only if
$p_1\Phi=0$ and the last condition means exactly
(\ref{m0onshuffle}).

\begin{proposition} Let
$ \{f_i\}: (M,\{m_i\})\to (M',\{m'_i\}) $ be an
$A_{\infty}$-algebra morphism of $C_{\infty}$-algebras. Then the
induced map of bar constructions $ \tilde{B}\{f_i\}$ is a map of
dg bialgebras if and only if each $f_i$ disappears on shuffles,
that is $\{f_i\}$ is a morphism of $C_{\infty}$-algebras.
\end{proposition}
\noindent {\bf Proof.} The map
$\Psi=\tilde{B}\{f_i\}\mu_{sh}-\mu_{sh}(\tilde{B}\{f_i\}\otimes
\tilde{B}\{f_i\}) $ is a coderivation. Thus, according to
universal property of $T^c(s^{-1}M)$ the map $\Psi$ is trivial if
and only if $p_1\Psi=0$ and the last condition means exactly
(\ref{f0onshuffle}).

Thus the bar functor maps the subcategory of $C_{\infty}$-algebras
to the category of dg bialgebras.

\subsection{Adjunctions}

The bar and cobar functors
$$
B:DGAlg\to DGCoalg,\  \Omega :DGCoalg\to DGAlg
$$
are adjoint and there exist standard weak equivalences $\Omega
B(A)\to A,\ \ C\to B\Omega C$. So $\Omega B(A)\to A$ is a {\it
free resolution} of a dga $A$.

If $A$ is commutative, the cobar-bar resolution is out of
category: $\Omega B(A)$ is not commutative.

In this case instead the cobar-bar functors we must use the
adjoint functors $\Gamma,\ {\cal A}$, see \cite{SS93}, which we
describe now.

For a commutative dg algebra the bar construction is a dg
bialgebra, so the restriction of the bar construction is the
functor $B:CDGAlg \to DGBialg$. Furthermore, the functor of
indecomposables $Q:DGBialg\to DGLieCoalg$ maps the category of dg
bialgebras to the category of dg Lie calgebras. Let $\Gamma$ be
the composition
$$
\Gamma: CDGAlg \stackrel{B}{\to}
DGBialg\stackrel{Q}{\to}DGLieCoalg.
$$
There is the adjoint of $\Gamma$ $ {\cal A}:DGLieCoalg \to CDGAlg,
$ which is dual to Chevalle-Eilenberg functor. There is the
standard weak equivalence ${\cal A}\Gamma A\to A$.

\subsection{Minimality}
\label{minimal}

Let $ \{f_i\}: (M,\{m_i\})\to (M',\{m'_i\})$ be a morphism of
$A_{\infty}$-algebras. It follows from (\ref{morphism}) that the
first component $f_1:(M,m_1)\to (M',m'_1)$ is a chain map.

A weak equivalence of $A_{\infty} $-algebras is defined as a
morphism $\{f_i\}$ for which $B(\{f_i\})$ is a weak equivalence of
dg coalgebras. The standard spectral sequence argument allows to
prove the following
\begin{proposition}
\label{mapweak}
 A morphism of $A_{\infty}$-algebras is a weak
equivalence if and only if it's first component $f_1:(M,m_1)\to
(M',m'_1)$ is a weak equivalence of chain complexes.
\end{proposition}
\begin{proposition}
\label{mapiso}
 A morphism of $A_{\infty}$-algebras is an
isomorphism if and only if it's first component $f_1:(M,m_1)\to
(M',m'_1)$ is an isomorphism.
\end{proposition}
\noindent {\bf Proof.} The components of opposite morphism $
\{g_i\}: (M',\{m'_i\})\to (M,\{m_i\})$ can be solved inductively
from the equation $\{g_i\}\{f_i\}=\{id_M,0,0,...\}. $
\begin{definition}
An $A_{\infty}$-algebra $(M,\{m_i\})$ we call minimal if $m_1=0$.
\end{definition}
In this case $(M,m_2)$ is \emph{ strictly} associative graded
algebra.

From the above propositions easily follows
\begin{proposition}
\label{aminimaliso}
 Each weak equivalence of minimal $A_{\infty}$-algebras is an isomorphism.
\end{proposition}

It is clear that all above is true for $C_{\infty}$-algebras, thus
\begin{proposition}
\label{cminimaliso}
 Each weak equivalence of minimal $C_{\infty}$-algebras is an isomorphism.
\end{proposition}

\begin{definition}
A minimal $A_{\infty}$-algebra ($C_{\infty}$-algebra)
$(M,\{m_i\})$ we call degenerate if it is isomorphic in the
category of $A_\infty$ $(C_\infty$) algebras
 to the graded (commutative) algebra $(M,m_2)$.
\end{definition}

\section{Minimal $A_{\infty}$  and $C_{\infty}$ algebras and
Hochschild and Harrison Cohomology} \label{hochschild}

Here we present the connection of the notion of minimal $A_\infty$
(resp. $C_\infty$)-algebra with Hochschild (resp. Harrison)
cochain complexes, studied in \cite{Kad88}, see also \cite{Kad07}.

Let $H$ be a graded algebra. Consider the Hochshild cochain
complex $C^{*,*}(H,H)$ which is bigraded in this case:
$$C^{n,m}(H,H)=Hom^m(H^{\otimes n},H),$$
where $Hom^m$ means homomorphisms of degree $m$.

 This bigraded complex
carries a structure of homotopy Gerstenhaber algebra, see
\cite{Kad88}, \cite{GJ94}, \cite{GV95}, \cite{Kad07}, which
consists of following structure maps:

(i) The Hochschild differential $\delta :C^{n-1,m}(H,H)\to
C^{n,m}(H,H)$ given by
$$
\begin{array}{ll}
\delta f(a_1\otimes ...\otimes a_n)=&a_1\cdot f(a_2\otimes ...\otimes a_n) \\
&+\sum_{k}\pm f(a_1\otimes ...\otimes a_{k-1}\otimes a_k\cdot
a_{k+1}\otimes ..\otimes a_n)\\
&\pm f(a_1\otimes ...\otimes a_{n-1})\cdot a_n;
\end{array}
$$
(ii) The $\smile $ product defined by
$$
f\smile g(a_1\otimes ...\otimes a_{n+m})=f(a_1\otimes ...\otimes
a_n)\cdot g(a_{n+1}\otimes ...\otimes a_{n+m}).
$$
(iii) The {\it brace} operations $f\{g_1,...,g_i\}$ which we write
as $f\{g_1,...,g_i\}=E_{1,i}(f;g_1,...,g_i)$,
$$ E_{1,i}:C^{n,m}\otimes C^{n_1,m_1}\otimes...\otimes
C^{n_i,m_i}\to C^{n+\sum n_t-i,m+\sum m_t},
$$
given by
\begin{equation}
\label{cup1gg}
\begin{array}{ll}
E_{1,i}(f;g_1,...,g_i)(a_1\otimes ...\otimes
a_{n+n_1+...+n_i-i})&\\
= \sum_{k_1,...,k_i} \pm f(a_1\otimes ...\otimes a_{k_1}\otimes
g_1(a_{k_1+1}\otimes ... \otimes a_{k_1+n_1})
\otimes &\\
\ \ \ ... \otimes a_{k_2}\otimes g_2(a_{k_2+1}\otimes ... \otimes
a_{k_2+n_2})\otimes a_{k_2+n_2+1}
\otimes ...\\
\ \ \ \otimes a_{k_i}\otimes g_i(a_{k_i+1}\otimes ... \otimes
a_{k_i+n_i})\otimes...\otimes a_{n+n_1+...+n_i-i}).
\end{array}
\end{equation}
The first brace operation $E_{1,1}$ has the properties of
Steenrod's $\smile_1$ product, so we use the notation
$E_{1,1}(f,g)=f\smile_1 g$. In fact this is Gerstenhaber's $f\circ
g$ product \cite{Ger63}, \cite{Ger64}.

Now let $(H,\{m_i\}$ be a minimal $A_{\infty}$-algebra, so
$(H,m_2)$ is an associative graded algebra with multiplication
$a\cdot b=m_2(a\otimes b)$.

Each operation $m_i$ can be considered as a Hochschild cochain
$m_i\in C^{i,2-i}(H,H)$. Let $m=m_3+m_4+...\in C^{*,2-*}(H,H)$.
The defining condition of $A_{\infty}$-algebra (\ref{ainf}) means
exactly $\delta m=m\smile_1 m$. So a minimal $A_{\infty}$-algebra
structure on $H$ in fact is a {\it twisting cochain} in the
Hochschild complex with respect to the $\smile_1$ product.

There is the notion of equivalence of such twisting cochains:
$m\sim m'$ if there exists
$p=p^{2,-1}+p^{3,-2}+...+p^{i,1-i}+...,\ ,\ p^{i,1-i}\in
C^{i,1-i}(H,H)$ such that
\begin{equation}
\label{eqv}
\begin{array}{l}

m-m'= \delta p+ p\smile p+p\smile_1m+ \\
 m'\smile_1p+
E_{1,2}(m';p,p)+ E_{1,3}(m';p,p,p)+ ...\ .
\end{array}
\end{equation}
\begin{proposition} Twisting cochains $m, m'\in C^{*,2-*}(H,H)$ are
equivalent if and only if $(H,\{m_i\})$ and $(H',\{m'_i\})$ are
isomorphic $A_{\infty}$-algebras.
\end{proposition}
\noindent{\bf Proof.} Indeed,
$$
\{p_i\}:(H,\{m_i\})\to (H,\{m'_i\})
$$
with $p_1=id,\ p_i=p^{i,1-i}$ is the needed isomorphism: the
condition (\ref{eqv}) coincides with the defining condition
(\ref{morphism}) of a morphism of $A_\infty$-algebras and the
Proposition \ref{mapiso} implies that this morphism is an
isomorphism.

This gives the possibility of perturbation of twisting cochain
without changing their equivalence class:
\begin{proposition}
Let $m$ be a twisting cochain (i.e. a minimal $A_{\infty}$-algebra
structure on $H$) and $ p\in C^{n,1-n}(H,H)$ be an arbitrary
cochain, then there exists a twisting cochain $\bar{m}$,
equivalent to $m$, such that $m_i=\bar{m}_i$ for $i\leq n$ and $
\bar{m}_{n+1}=m_{n+1}+\delta p$.
\end{proposition}
\noindent{\bf Proof.} The twisting cochain $\bar{m}$ can be solved
inductively from the equation (\ref{eqv}).

\begin{theorem}\label{hoch3} Suppose for a graded algebra $H$ Hochschild cohomology
$Hoch^{n,2-n}(H,H)=0$ for $n\geq 3$. Then each $m\sim 0$, that is
each minimal $A_{\infty}$-algebra structure on $H$ is degenerate.
\end{theorem}

\noindent {\bf Proof.} From the equality $\delta m=m\smile_1 m$ in
dimension $4$ we obtain $\delta m_3=0$ that is $m_3$ is a cocycle.
Since $Hoch^{3,-1}(H,H)=0$ there exists $p^{2,-1}$ such that
$m_3=\delta p^{2,-1}$. Perturbing our twisting cochain $m$ by
$p^{2,-1}$ we we obtain new twisting cochain
$\bar{m}=\bar{m}_3+\bar{m}_4+...$ equivalent to $m$ and with
$\bar{m}_3=0$. Now the component $\bar{m}_4$ becomes a cocycle,
which can be killed using $Hoch^{4,-2}(H,H)=0$ etc.


Suppose now $(H,\mu)$ is a commutative graded algebra. The
Harrison cochain complex $\bar{C}^*(H,H)$ is defined as a
subcomplex of the Hochschild complex consisting of cochains which
disappear on shuffles. If $(H,\{m_i\})$ is a $C_{\infty}$-algebra
then the twisting element $m=m_3+m_4+...$ belongs to Harrison
subcomplex $\bar{C}^*(H,H)\subset C^*(H,H)$ and we have the
\begin{theorem}\label{harr3} Suppose for a graded commutative algebra $H$
Harrison cohomology $Harr^{n,2-n}(H,H)=0$ for $n\geq 3$. Then each
$m\sim 0$, that is each minimal $C_{\infty}$-algebra structure on
$H$ is degenerate.
\end{theorem}


\section{$A_{\infty}$-algebra structure in homology}


Let $(A,d,\mu )$ be a dg algebra and $(H(A),\mu ^{*})$ be it's
homology algebra. Although the product in $H(A)$ is associative,
there appears a structure of a (generally nondegenerate) minimal
$A_{\infty}$-algebra, which can be considered as an $A_{\infty}$
\emph{deformation} of $(H(A),\mu^*)$, \cite{Kad07} . Namely, in
\cite{Kad76}, \cite {Kad80} the following result was proved (see
also \cite{Smi80}, \cite{GS86}):

\begin{theorem}
Suppose for a dg algebra $A$ all homology modules $H^i(A)$ are
free.

Then there exist: a structure of minimal $A_{\infty}$-algebra
$(H(A),\{m_i\})$ on $H(A)$ and a weak equivalence of
$A_{\infty}$-algebras
$$
\{f_i\}:(H(A),\{m_i\})\rightarrow (A,\{d,\mu ,0,0,...\})
$$
such, that $m_1=0$, $m_2=\mu ^{*}$, $f_1^{*}=id_{H(A)}$.

Furthermore, for a dga map $f:A\to A'$ there exists a morphism of
$A_{\infty}$-algebras $\{f_i\}:(H(A)\{m_i\})\to (H(A')\{m'_i\})$
with $f_1=f^*$.
\end{theorem}

 Such a
structure is unique up to isomorphism in the category of
$A_{\infty}$-algebras: if $(H(A),\{m_i\})$ and $(H(A),\{m'_i\})$
are two such $A_{\infty}$-algebra structures on $H(A)$ then for
$id:A\to A$ there exists $\{f_i\}:(H(A)\{m_i\})\to (H(A)\{m'_i\})$
with $f_1=id$, so, since of Proposition \ref{mapiso} $\{f_i\}$ is
an isomorphism.

Let us look at the first new operation $m_3:H(A)\otimes
H(A)\otimes H(A)\to H(C)$. Let $f_1:H(A)\to A$ be a cycle-choosing
homomorphism: $f_1(a)\in a\in H(A)$. This map is not
multiplicative but $ f_1(a\cdot b)-f_1(a)\cdot f(b)\sim 0\in C $
so there exists $f_2:H(A)\otimes H(A)\to A$ s.t. $f_1(a\cdot
b)-f_1(a)\cdot f(b)=\partial f_2(a\otimes b)$. We define
$m_3(a\otimes b\otimes c)\in H(A)$ as the homology class of the
cycle
$$
f_1(a)\cdot f_2(b\otimes c)\pm f_2(a\cdot b\otimes c)\pm
f_2(a\otimes b\cdot c)\pm f_2(a\otimes b)\cdot f_1(c).
$$
From this description immediately follows the connection of $m_3$
with Massey product: If $a,b,c\in H(A)$  is a  Massey triple, i.e.
if $a\cdot b=b\cdot c=0$, then $m_3(a\otimes b\otimes c)$ belongs
to the Massey product $<a,b,c>$. This gives examples of gd
algebras with essentially nontrivial homology $A_\infty$-algebras.

 \subsection{Main examples  and applications}

Taking $A=C^*(X)$, the cochain dg algebra of a 1-connected space
$X$, we obtain an $A_{\infty}$-algebra structure
$(H^*(X),\{m_i\})$ on cohomology algebra $H^*(X)$.

Cohomology algebra equipped with this additional structure carries
more information then just the cohomology algebra. Some
applications of this structure are given in \cite{Kad80} ,
\cite{Kad93}. For example the cohomology $A_{\infty}$-algebra
$(H^*(X),\{m_i\})$ determines cohomology of the loop space
$H^*(\Omega X)$ when just the algebra $(H^*(X),m_2)$ does not:

\begin{theorem} $H(\tilde{B}(H^*(X),\{m_i\}))=H^*(\Omega X)$.
\end{theorem}

Taking $A=C_*(G)$, the chain dg algebra of a topological group
$G$, we obtain an $A_{\infty}$-algebra structure
$(H_*(G),\{m_i\})$ on the Pontriagin algebra $H_*(G)$.
 The homology
$A_{\infty}$-algebra $(H_*(G),\{m_i\})$ determines homology of the
classifying space $H_*(B_G)$ when just the Pontriagin algebra
$(H_*(G),m_2)$ does not:

\begin{theorem} $H(\tilde{B(H_*(G)},\{m_i\}))=H_*(B_G)$.
\end{theorem}



  \section{$C_{\infty}$-algebra structure in homology
  of a commutative dg algebra}

There is a commutative version of the above main theorem,
see\cite{Kad88h}, \cite{Kad93}, \cite{Markl}:
\begin{theorem}\label{cinf}
Suppose for a commutative dg algebra $A$ all homology $R$-modules
$H^i(A)$ are free.

Then there exist: a structure of minimal $C_{\infty}$-algebra
$(H(A),\{m_i\})$ on $H(A)$ and a weak equivalence of
$C_{\infty}$-algebras
$$
\{f_i\}:(H(A),\{m_i\})\rightarrow (A,\{d,\mu ,0,0,...\})
$$
such, that $m_1=0$, $m_2=\mu ^{*}$, $f_1^{*}=id_{H(A)}$.

Furthermore, for a cdga map $f:A\to A'$ there exists a morphism of
$C_{\infty}$-algebras $\{f_i\}:(H(A)\{m_i\})\to (H(A')\{m'_i\})$
with $f_1=f^*$.
\end{theorem}
 Such a
structure is unique up to isomorphism in the category of
$C_{\infty}$-algebras.

Bellow we present some applications of this $C_\infty$-algebra
structure in rational homotopy theory.



 \section{Applications in Rational Homotopy Theory}

\subsection{Classification of rational homotopy types}

Let $X$ be a 1-connected space. In the case of rational
coefficients there exist Sullivan's {\it commutative} cochain
complex $A(X)$ of $X$. It is well known that the weak equivalence
type of cdg algebra $A(X)$ determines the rational homotopy type
of $X$: 1-connected $X$ and $Y$ are {\it rationally homotopy
equivalent} if and only if $A(X)$ and $A(Y)$ are weekly homotopy
equivalent cdg algebras. Indeed, in this case $A(X)$ and $A(Y)$
have {\it isomorphic} minimal models $M_X\approx M_Y$, and this
implies that $X$ and $Y$ are rationally homotopy equivalent. This
is the key geometrical result of Sullivan which we are going to
exploit bellow.

Now we take $A=A(X)$ and apply the Theorem \ref{cinf}. Then we
obtain on $H(A)=H^*(X,Q)$ a structure of minimal $C_{\infty}$
algebra $(H^*(X,Q),\{m_i\})$ which we call {\it rational
cohomology $C_{\infty}$-algebra of $X$}.

Generally isomorphism of rational cohomology algebras $H^*(X,Q)$
and $ H^*(Y,Q)$ does not imply homotopy equivalence $X\sim Y$ even
rationally. We claim that $(H^*(X,Q),\{m_i\})$ is {\it complete}
rational homotopy invariant:
\begin{theorem} 1-connected $X$ and $X'$ are rationally
homotopy equivalent if and only if
 $(H^*(X,Q),\{m_i\})$ and $(H^*(X',Q),\{m'_i\})$ are
isomorphic as $C_{\infty}$-algebras.
\end{theorem}
\noindent{\bf Proof.} Suppose $X\sim X'$, then $A(X)$ and $A(X')$
are weak equivalent, that is there exists a cgda $A$ and weak
equivalences $A(X)\leftarrow A\to A(X')$. This implies weak
equivalences of corresponding homology $C_{\infty}$-algebras
$$
(H^*(X,Q),\{m_i\})\leftarrow (H^*(A),\{m_i\})\to
(H^*(X',Q),\{m'_i\}),
$$
which since of minimality both are isomorphisms.

Conversely, suppose $(H^*(X,Q),\{m_i\})\approx
(H^*(X',Q),\{m'_i\})$. Then  $${\cal A}
QB(H^*(X,Q),\{m_i\})\approx {\cal A} QB(H^*(X',Q),\{m'_i\}).$$
Denote this cdga as $A$. Then we have weak equivalences of CGD
algebras
$$
A(X)\leftarrow {\cal A} \Gamma A(X)\leftarrow A\to {\cal A} \Gamma
A(X') \to A(X').
$$

This theorem in fact classifies rational homotopy types with given
cohomology algebra $H$ as all possible minimal
$C_{\infty}$-algebra structures on $H$ modulo $C_{\infty}$
isomorphisms.

\noindent{\bf Example.} Here we describe an example which we will
use to illustrate the results of this and forthcoming sections.

We consider the following commutative graded algebra. It's
underline graded $Q$-vector space has the generators: generator
$e$ of dimension 0, generators $x,\ y$ of dimension 2, and
generator $z$ of dimension 5, so
\begin{equation}\label{example}
H^*=\{H^0=Q_e,\ 0,\ H^2=Q_x\oplus Q_y,\ 0,\ 0,\ H^5=Q_z,\ 0,\ 0,\
...\ \},
\end{equation}
and the multiplication is trivial by dimensional reasons, with
unit $e$. In fact
$$
H^*=H^*(S^2\vee S^2\vee S^5,Q).
$$
This example was considered in \cite{HS} and there was shown that
there are just two rational homotopy types with such cohomology
algebra.

The same result can be obtained from our classification.

What minimal $C_\infty$-algebra structures are possible on $H^*$?

By dimensional reasons only one nontrivial operation $
m_3:H^2\otimes H^2\otimes H^2\to H^5 $ is possible.

The specific condition of $C_\infty$-algebra, namely the
disappearance on shuffles implies that
$$
m_3(x,x,x)=0,\ m_3(y,y,y)=0,\ m_3(x,y,x)=0,\ m_3(y,x,y)=0
$$
and
$$
m_3(x,x,y)=m_3(y,x,x),\ \ m_3(x,y,y)=m_3(y,y,x).
$$
Thus each $C_\infty$-algebra structure on $H^*$ is characterized
by a couple rational numbers $p,\ q$,
$$
m_3(x,x,y)=pz,\ \ m_3(x,y,y)=qz.
$$
So let us write an arbitrary minimal $C_\infty$-algebra structure
on $H^*$ as a column vector
$\left(\begin{array}{c}p\\q\end{array}\right)$.

Now let us look at the structure of an isomorphism of
$C_\infty$-algebras
$$
\{f_i\}:(H^*,m_3)\to (H^*,m'_3).
$$
Again by dimensional reasons just one component $f_1:H^*\to H^* $
is possible, which in it's turn consists of two isomorphisms
$$
f_1^2:H^2=Q_x\oplus Q_y\to H^2=Q_x\oplus Q_y,\ \ f_1^5:H^5=Q_z\to
H^5=Q_z.
$$
The first one is represented by a nondegenerate matrix
$A=\left(\begin{array}{cc}a&c\\b&d\end{array}\right)$,
$$
f^2_1(x)=ax\oplus by,\ \ f^2_1(y)=cx\oplus dy,
$$
and the second one by a nonzero rational number $r$,
$f^5_2(z)=rz$.

Calculation shows that the condition $f^5_1m_3=m'_3(f^2_1\otimes
f^2_1\otimes f^2_1)$, to which degenerates the defining condition
of an $A_\infty$-algebra morphism (\ref{morphism}) looks as
$$
r\left(\begin{array}{c}p\\q\end{array}\right)=det\
A\left(\begin{array}{cc}a&b\\c&d\end{array}\right)
\left(\begin{array}{c}p'\\q'\end{array}\right).
$$
This condition shows that two minimal $C_\infty$-algebra
structures $m_3=\left(\begin{array}{c}p\\q\end{array}\right)$ and
$m'_3=\left(\begin{array}{c}p'\\q'\end{array}\right)$ are
isomorphic if and only if they are tied with nondegenerate linear
transformation.

Thus that there exist just two isomorphism classes of minimal
$C_\infty$-algebra on $H^*$: the trivial one $(H^*,m_3=0)$ and the
nontrivial one $(H^*,m_3\neq 0)$. So we have just two rational
homotopy types whose rational cohomology is $H^*$. We denote them
$X$ and $Y$ respectively and analyze in next sections.

\vspace{5mm}

 Below we give some applications of cohomology $C_\infty$-algebra in various
 problems of rational homotopy theory.

\subsection{Formality}

Among rational homotopy types with given cohomology algebra, there
is one called {\it formal} which is "formal consequence of it's
cohomology algebra" (Sullivan). Explicitly this is the type whose
minimal model $M_X$ is isomorphic to the minimal model of
cohomology $H^*(X,Q)$.

Our $C_{\infty}$ model implies the following criterion of
formality:

\begin{theorem}\label{formal} X is formal if and only if its cohomology
$C_{\infty}$-algebra is degenerate, i.e. it is $C_{\infty}$
isomorphic one with $m_{\geq 3}=0$.
\end{theorem}

Bellow we deduce using this criterion some known results about
formality.

\noindent 1. A commutative graded 1-connected algebra $H$ is
called {\it intrinsically formal} if there is only one homotopy
with cohomology algebra $H$, of course the formal one.

The above Theorem \ref{harr3} immediately implies the following
sufficient condition for formality due to Tanre \cite{Tan85}:

\begin{theorem}  If for a 1-connected graded $Q$-algebra $H$ one has
$$Harr^{k,k-2}(H,H)=0,\ k=3,4,...$$ then
$H$ is intrinsically formal, that is there exists only one
rational homotopy type with $H^*(X,Q)\approx H$.
\end{theorem}

\noindent 2. The following theorem of Halperin and Stasheff from
\cite{HS} is an immediate result of our criterion:
\begin{theorem}\label{halsta} A commutative graded $Q$-algebra of
type
$$
H=\{H^0=Q,0,0,...,0,H^n,H^{n+1},...,H^{3n-2},0,0,...\}
$$
is intrinsically formal
\end{theorem}
\noindent{\bf Proof.} Since $deg\ m_i=2-i$ there is no room for
operations $m_{i>2}$, indeed the shortest range is $m_3:H^n\otimes
H^n\otimes H^n\to H^{3n-1}=0. $

\vspace{5mm}

\noindent 3. From the Theorem \ref{formal} easily follows the
\begin{theorem} Any  1-connected commutative graded algebra $H$ with
$H^{2k}=0$ is intrinsically formal.
\end{theorem}
\noindent {\bf Proof.} Any $A_{\infty}$-operation $m_i$ has degree
$2-i$, thus
$$
m_i:H^{2k_1+1}\otimes ...\otimes H^{2k_i+1}\to
H^{2(k_1+...+k_i+1)}=0.
$$
Thus any $C_{\infty}$ operation is trivial too.

From this follows one result of Baues: any space whose even
dimensional cohomologiies are trivial has rational homotopy type
of wedge of spheres. Indeed, such algebra is realized as a wedge
of spheres and since of intrinsical formality this is the only
homotopy type.


\noindent{\bf Example.} The algebra $H^*$ from the example of
previous section is not intrinsically formal since there are two
homotopy types, $X$ and $Y$, with $H^*(X,Q)=H^*=H^*(Y)$. The space
$X$  is formal (and actually $X=S^2\vee S^2\vee S^5$), since it's
cohomology $C_\infty$-algebra $(H^*,m_3=0)$ is trivial. But the
space $Y$ is not: it's cohomology $C_\infty$-algebra $(H^*,m_3\neq
0)$ is not degenerate.

We remark here that the formal type is represented by $X=S^2\vee
S^2\vee S^5$ and it is possible to show that the nonformal one is
represented by $Y=S^2\vee S^2\cup_{f:S^4\to S^2\vee S^2} e^5$,
where the attaching map $f$ is a nontrivial element from
$\pi_4(S^2\vee S^2)\otimes Q$.


\subsection{Rational homotopy groups}

Since the cohomology $C_{\infty}$-algebra $(H^*(X,Q),\{m_i\})$
determines the rational homotopy type it must determine the
rational homotopy groups $\pi_i(X)\otimes Q$ too. We present a
chain complex whose homology is $\pi_i(X)\otimes Q$. Moreover the
Lie algebra structure is determined as well.

For cohomology $C_{\infty}$-algebra $(H^*(X,Q),\{m_i\})$ the bar
construction $B(H^*(X,Q),\{m_i\})$ is dg bialgebra. Acting on this
bialgebra by the functor $Q$ of indecomposables we obtain a dg Lie
coalgebra.

On the other hand rational homotopy groups $\pi_*(\Omega X)\otimes
Q$ form a graded Lie algebra with respect to Whiethead product.
Thus it's dual cohomotopy groups $\pi^*(\Omega X,Q)=(\pi_*(\Omega
X)\otimes Q)^*$ form a graded Lie coalgebra.

\begin{theorem} Homology of dg Lie coalgebra $QB(H^*(X,Q),\{m_i\})$ is isomorphic to
cohomotopy Lie coalgebra $\pi^*(\Omega X,Q)$.
\end{theorem}
\noindent{\bf Proof.} The theorem follows from the sequence of
graded Lie coalgebra isomorphisms:
$$
\begin{array}{ll}
\pi^*(\Omega X,Q)\approx  &(\pi_*(\Omega X,Q))^*\approx
 (PH_*(\Omega
x,Q)^* \approx QH^*(\Omega X,Q)\approx \\
&QH(B(A(X)) \approx QH(\tilde{B}(H^*(X,Q),\{m_i\}) \approx\\
 &H(Q\tilde{B}(H^*(X,Q),\{m_i\}).
\end{array}
$$

\noindent{\bf Example.} For the algebra $H^*$ from the previous
examples the complex $QB(H^*)$ in low dimensions looks as
$$
0\to Q_x\oplus Q_y\stackrel{0}{\to}Q_{x\otimes x}\oplus
Q_{x\otimes y}\oplus Q_{y\otimes y}\stackrel{0}{\to} Q_{x\otimes
x\otimes y}\oplus Q_{x\otimes y\otimes y}\stackrel{d=m_3}{\to}
Q_z\oplus ... \ .
$$
The differential $d=m_3$ is trivial for the formal space $X$ and
is nontrivial for $Y$. Thus for both rational homotopy types we
have
$$
\pi^2=H^1(QB(H^*))=2Q,\ \ \pi^3=H^2(QB(H^*))=3Q,
$$
and
$$
\begin{array}{l}
\pi^4(X)=H^3(QB(H^*),d=0)=2Q,\\
\pi^4(Y)=H^3(QB(H^*),d\neq 0)=Ker\
d=Q.
\end{array}
$$


\subsection{Realization of homomorphisms}

Let $G:H^*(X,Q)\to H^*(Y,Q)$ be a homomorphism of cohomology
algebras. When this homomorphism is realizable as a map of
rationalizations $g:Y_Q\to X_Q ,\ \ f^*=F$? In the case when $G$
is an isomorphism this question was considered in \cite{HS}. It
was considered also in \cite{Vig81}. The following theorem gives
the complete answer:

\begin{theorem}\label{realizat} A homomorphism $G$ is realizable if and only if it is extendable to a
$C_{\infty}$-map
$$
\{g_1=G,g_2,g_3,...\}:(H^*(X,Q),\{m_i\})\to
(H^*(Y,Q),\{m'_i\}).
$$
\end{theorem}
\noindent{\bf Proof.} One side of is consequence of the last part
of Theorem \ref{cinf}.

To show the other side we use Sullivan's minimal models $M_X$ and
$M_Y$ of $A(X)$ and $A(Y)$. It is enough to show that the
existence of $\{g_i\}$ implies the existence of cdg algebra map
$g:M_Y\to M_X$.

So we have $C_{\infty}$-algebra maps
$$
M_X\stackrel{\{f_i\}}{\leftarrow}
(H^*(X,Q),\{m_i\})\stackrel{\{g_i\}}{\to}
(H^*(y,Q),\{m'_i\})\stackrel{\{f'_i\}}{\to}M_Y.
$$
Recall the following property of a minimal cdg algebra $M$: for a
weak equivalence of cdg algebras $\phi:A\to B$ and a cdg algebra
map $f:M\to B$ there exists a cdg algebra map $F:M\to A$ such that
$\phi F$ is homotopic to $f$. Using this property it is easy to
show the existence of a cdga map $\beta:M_X\to {\cal A}QB(M_X)$,
the right inverse of the standard map $\alpha:{\cal A}QB(M_X)\to
M$. Composing this map with ${\cal A}QB({\{f'_i\}}){\cal
A}QB({\{g_i\}})$ we obtain a cdga map
$$
{\cal A}QB({\{f'_i\}}){\cal A}QB({\{g_i\}})\beta:M_X\to M_Y.
$$
From this theorem immediately follows the {\begin{corollary} For
formal $X$ and $Y$ each $G:H^*(X,Q)\to H^*(Y,Q)$ is realizable.
\end{corollary}
\noindent{\bf Proof.} In this case $\{G,0,0,...\}$ is a
$C_{\infty}$ extension f $G$.

\noindent{\bf Example.} Consider the homomorphism
$$
G:H^*(X)=H^*(Y)\to H^*(S^5)
$$
induced by the standard imbedding $g:S^5\to X=S^2\vee S^2\vee
S^5$. Of course $G$ is realizable as $g:S^5\to X$ but not as
$S^5\to Y$. Indeed, for such realizability, according to Theorem
\ref{realizat}, we need a $C_\infty$-algebra morphism
$$
\{g_i\}:(H^*,\{0,0,m_3,0,...\}\to (H^5(S^5,Q),\{0,0,0,...\})
$$
with $g_1=G$. By dimensional reasons all the components
$g_2,g_3,...$ all are trivial, so this morphism looks as
$\{G,0,0,...\}$. But this collection is not a morphism of
$C_\infty$-algebras since the condition $Gm_3=0$, to which
degenerates the defining condition (\ref{morphism}) of an
$A_\infty$-algebra morphism, is not satisfied.

A. Razmadze Mathematical Institute,

1, M. Alexidze Str., Tbilisi, 0193, Georgia

kade@rmi.acnet.ge

\end{document}